\documentclass[12pt]{article}

\usepackage{latexsym,amssymb,amsmath,amsthm,mathrsfs,euscript,graphicx}
\usepackage{arydshln}
\usepackage{hyperref}
\usepackage{enumitem}

\newtheorem{theorem}{Theorem}[section]
\newtheorem{lemma}[theorem]{Lemma}
\newtheorem{corollary}[theorem]{Corollary}

\newcommand\numberthis{\addtocounter{equation}{1}\tag{\theequation}}

\def\barroman#1{\sbox0{#1}\dimen0=\dimexpr\wd0+1pt\relax
  \makebox[\dimen0]{\rlap{\vrule width\dimen0 height 0.06ex depth 0.06ex}%
    \rlap{\vrule width\dimen0 height\dimexpr\ht0+0.03ex\relax 
            depth\dimexpr-\ht0+0.09ex\relax}%
    \kern.5pt#1\kern.5pt}}
    
\begin{document}

\setlength{\baselineskip}{1.25 \baselineskip}

\def \today {\number \day \ \ifcase \month \or January\or February\or
  March\or April\or May\or June\or July\or August\or
  September\or October\or November\or December\fi\
  \number \year}


\title{\sffamily On random digraphs and cores}

\author{
   \textsc{Esmaeil Parsa\footnotemark~\footnotemark} \\[0.25em]
    {\small\textit{Department of Mathematics}} \\[-0.25em]
    {\small\textit{Islamic Azad University (Parand Branch) }} \\[-0.25em]
    {\small\textit{Parand New Town, Iran}} \\[-0.1em]
    {\small\texttt{esmaeil.parsa@yahoo.com}}
    \and\addtocounter{footnote}{-2}
  \textsc{P. Mark Kayll\footnotemark} \\[0.25em]
    {\small\textit{Department of Mathematical Sciences}} \\[-0.25em]
    {\small\textit{University of Montana}} \\[-0.25em]
    {\small\textit{Missoula MT 59812, USA}} \\[-0.1em]
    {\small\texttt{mark.kayll@umontana.edu}}
   }

\date{\small \today}

\maketitle

%
\renewcommand{\thefootnote}{}

%
%
\renewcommand{\thefootnote}{\fnsymbol{footnote}}

\addtocounter{footnote}{1}
\footnotetext{Partially supported by a grant from the Simons Foundation (\#279367 to Mark Kayll)} 
\addtocounter{footnote}{1}
\footnotetext{Partially supported by a 2017 University of Montana
  Graduate Student Summer Research Award funded by the George and Dorothy Bryan
Endowment}
\footnotetext{This work forms part of the author's PhD dissertation~\cite{Parsa}.}

  \begin{abstract}

\noindent
An acyclic homomorphism of a digraph $C$ to a digraph $D$ is a
function  $\rho\colon V(C)\to V(D)$ such that for every arc $uv$ of
$C$, either $\rho(u)=\rho(v)$, or $\rho(u)\rho(v)$ is an arc of $D$
and for every vertex $v\in V(D)$, the subdigraph of $C$ induced by
$\rho^{-1}(v)$ is acyclic. A digraph $D$ is a core if the only acyclic
homomorphisms of $D$ to itself are automorphisms. In this paper, we
prove that for certain choices of $p(n)$, random digraphs 
$D\in D(n,p(n))$  are asymptotically almost surely cores. For
digraphs, this mirrors a result from 
[A. Bonato and P. Pra\l{}at, The good, the bad, and the great:
homomorphisms and cores of random graphs, 
\textit{Discrete Math.}, \textbf{309} (2009), no.\,18, 
5535--5539; MR2567955]  concerning random graphs and cores.  


   \end{abstract}

\renewcommand{\thefootnote}{\arabic{footnote}}
\addtocounter{footnote}{-2}

\section{Introduction}
\label{S:1}
In this paper, we follow~\cite{digraph} and~\cite{graph} for
definitions and terminology. Our digraphs are simple, i.e., loopless and
without multiple arcs. However, we allow two vertices $u,v$ to be
joined by two oppositely directed arcs, $uv$ and $vu$. By a \textit{cycle}, we
always mean a directed cycle in the digraph case. For a natural number
$n$ and $0\leq p \leq 1$, a digraph $D\in D(n,p)$ is defined to be a
digraph on $n$ vertices (we use $V(D)=[n]=\{1,2,\dots,n\}$) where each
ordered pair of vertices is joined by an arc with probability $p$,
with the arcs chosen independently. Note that if $D$ is any particular
digraph on $n$ vertices, then the probability of obtaining $D$ is
$p^{|A(D)|}(1-p)^{n(n-1)-|A(D)|}$. 

If $\mathscr{Q}$ is any digraph property (e.g., contains a
$\overleftrightarrow{K}_3$, is connected, etc.), we say that $D\in
D(n,p(n))$ has property $\mathscr{Q}$  ($D\in \mathscr{Q}$) 
\textit{a.a.s.\ (asymptotically almost surely)} if $P(D\in \mathscr{Q})\rightarrow 1$
as $n\rightarrow \infty$. We use $v_C$ and $a_C$ to denote $|V(C)|$
and $|A(C)|$,  respectively, for a digraph $C$. We
sometimes use the asymptotic notations $a_n\ll b_n$ and $a_n\asymp b_n$ to
denote $a_n=o(b_n)$ and $a_n=\Theta(b_n)$, respectively, for positive
sequences $(a_n)$ and $(b_n)$. 

The \textit{maximum density} of $D$ is $m(D):=\max
\{\frac{a_C}{v_C}:C\;\text{is a subdigraph of}\; D \;\text{and}\;
v_C>0\}$. Let $\mathscr{Q}$ be a nontrivial digraph property (a
property that is not satisfied by all or no digraphs). We say that
$\mathscr{Q}$ is \textit{monotone increasing} if $D\in \mathscr{Q}$
implies that $C\in \mathscr{Q}$ for every digraph $C$ on the same set
of vertices containing $D$ as a subdigraph. Let $\mathscr{Q}$ be a
nontrivial monotone increasing digraph property, $(\hat{p}_n)$ a
sequence of probabilities, and $D\in D(n,p(n))$. Then $(\hat{p}_n)$ is
a \textit{threshold} for $\mathscr{Q}$ if  
   \begin{align*}
        P(D\in \mathscr{Q})\rightarrow \begin{cases}
        0 & \text{if}\;\; p(n)\ll \hat{p}_n\\
        1 & \text{if}\;\; p(n)\gg \hat{p}_n
        \end{cases}
    \end{align*}
as $n\rightarrow \infty$. 

The following assertion is a digraph analogue of \cite[Theorem
3.4]{Rgraph} and can be proved following the same technique.  

\begin{theorem} For an arbitrary digraph $C$ with at least one arc,
\begin{align*}
    \lim_{n\rightarrow \infty}P(C\subseteq D\in D(n,p(n)))=\begin{cases}
    0   &\quad \text{if}\;p(n)\ll n^{-1/m(C)}\\
    1 & \quad \text{if}\; p(n)\gg n^{-1/m(C)}.
    \end{cases}
\end{align*}
\label{Threshold}
\end{theorem}

\section{Asymptotic properties of random digraphs}
\label{S:2}
We begin with Chernoff's inequality, which is used
extensively in the proof of Lemma~\ref{Incomparable}. Here $X\in B(n,p)$
indicates that $X$ is a binomial random variable with parameters $n$
and $p$, with $n$ being the number of trials and $p$ the success
probability of each trial.

\begin{theorem}[Chernoff's inequality~\cite{Rgraph}]  If $X\in B(n,p)$ and $\lambda=np$, then, with $\rho(x)=(1+x)\log(1+x)-x$ for $x\geq -1$ (and $\rho(x)=\infty$ for $x<-1$), we have 
\begin{align*}
      P(X\geq E(X)+t)\leq  \exp(-\lambda\rho(t/\lambda))\leq \exp\Big(-\frac{t^2}{2(\lambda+t/3)}\Big)\;\;\text{for}\;\;t\geq 0,
      \end{align*}
      and
      \begin{align*}
    P(X\leq E(X)-t)\leq \exp(-\lambda\rho(-t/\lambda))\leq \exp\Big(-\frac{t^2}{2\lambda}\Big)\;\;\text{for}\;t\geq 0.
\end{align*}
\label{Chernoff's Inequality}
\end{theorem}
One immediate consequence of Theorem \ref{Chernoff's Inequality} is

\begin{corollary}[\cite{Rgraph}]
If $X\in B(n,p)$ and $\epsilon>0$, then
\begin{align*}
    P(|X-E(X)|\geq \epsilon E(X))\leq 2\exp(-\rho(\epsilon)E(X)).
\end{align*}
In particular, if $\epsilon\leq 3/2$, then
\begin{align*}
     P(|X-E(X)|\geq \epsilon E(X)\leq 2 \exp\Big(-\frac{\epsilon^{2}E(X)}{3}\Big).
\end{align*}
\label{Chernoffresult}
\end{corollary}
In order to prove the main result of this paper---Theorem~\ref{main-thm}---we need several lemmas, collected together in the following result. This extends Lemma~1 in~\cite{Bonato09} to random digraphs.
\begin{lemma} If $n^{-1/9}\log^2 n<p=p(n)<1-n^{-1/9}\log^2 n$, then a.a.s.\ $D\in D(n,p)$ has the following properties:
\begin{enumerate}[label=(\alph*),topsep=0pt]
    \item the number of neighbours of a vertex of $D$ is at least $n(2p-p^2)(1-o(1))$ and at most $n(2p-p^2)(1+o(1))$;
    \item every pair of distinct vertices of $D$ has at least $np^2(2-p)^2(1-o(1))$ and at most $np^2(2-p)^2(1+o(1))$ common neighbours;
    \item the largest acyclic subdigraph of $D$ has fewer than $n^{1/9}$ vertices;
    \item each set of $k$ vertices, where $k\geq k_0=k_0(n)=n^{1/9}\log^2n/2$, induces a subdigraph with at most $2p\binom{k}{2}(1+o(1))$ arcs;
    \item in each set of $k$ disjoint pairs of vertices $\{\{v_i,w_i\}\}$, for $i\in [k]$ where $k\geq k_1=k_1(n)=n^{1/9}\log^2n$,  there are at least $2(1-(1-p)^4)\binom{k}{2}(1+o(1))$ pairs $(i,j)$ such that at least one of $v_iv_j,v_iw_j,w_iv_j,w_iw_j$ is an arc of $D$.
\end{enumerate}
\label{Incomparable}
\end{lemma}
\begin{proof}

    (a)  Let $v$ be an arbitrary vertex of $D\in D(n,p)$. We define the random variable $X$ as $X=|N_D(v)|$. We have  
\begin{align*}
    E(X)=(n-1)[1-(1-p)^2]=(n-1)(2p-p^2)=n(2p-p^2)-O(1).
\end{align*}
Using Corollary~\ref{Chernoffresult} with $\epsilon=\log n/\sqrt{n(2p-p^2)}$ we have
\begin{align*}
    P(X\geq n(2p-p^2)+\sqrt{n(2p-p^2)}\log n\; \text{or}\; X &\leq n(2p-p^2)-\sqrt{n(2p-p^2)}\log n)\\ &
    \leq 2 \exp{(-\frac{\log^2 n}{3})}.
\end{align*}

Now, suppose that the random variable $Y$ counts all the vertices having at least $[ n(2p-p^2)+\sqrt{n(2p-p^2)}\log n]$ or at most $[n(2p-p^2)-\sqrt{n(2p-p^2)}\log n]$ neighbours. Using Markov's inequality, we have
\begin{align*}
    P(Y=0)=1-P(Y\geq 1)\geq 1-E(Y)\geq 1-2n\exp{(-\frac{\log^2
        n}{3})}\rightarrow 1\;\;\text{as } n\rightarrow\infty.
\end{align*}
So a.a.s.\ the number of neighbours of every vertex of $D\in D(n,p)$ lies between $n(2p-p^2)(1-o(1))$ and $n(2p-p^2)(1+o(1))$.

(b) Let $v_1$ and $v_2$ be two distinct vertices of $D\in D(n,p)$ and let $X$ count their common neighbours. Then
\begin{align*}
    E(X)=(n-2)[1-(1-p)^2][1-(1-p)^2]=(n-2)p^2(2-p)^2=np^2(2-p)^2-O(1).
\end{align*}
Using Corollary~\ref{Chernoffresult} with $\epsilon=\log n/\sqrt{np^2(2-p)^2}$, we have
\begin{align*}
   P(X\geq np^2(2-p)^2+\sqrt{np^2(2-p)^2}\log n \; & \text{or}\; X\leq np^2(2-p)^2-\sqrt{np^2(2-p)^2)}\log n) & \\ \leq 2\exp{(-\frac{\log^2 n}{3})}.
\end{align*}

Now, suppose that $Y$ counts all pairs of vertices having at least $[ np^2(2-p)^2+\sqrt{np^2(2-p)^2}\log n]$ or at most $[np^2(2-p)^2-\sqrt{np^2(2-p)^2}\log n]$ common neighbours. Then
\begin{align*}
    P(Y=0)=1-P(Y\geq 1) & \geq 1-E(Y)\geq 1-\binom{n}{2}2\exp{(-\frac{\log^2 n}{3})} \\ & 
    =1-O(n^2)\exp{(-\frac{\log^2 n}{3})}\rightarrow 1\;\;\text{as } n\rightarrow\infty.
\end{align*}
So a.a.s.\ the number of common neighbours of any two distinct vertices lies  between $np^2(2-p)^2(1-o(1))$ and $np^2(2-p)^2(1+o(1))$.

(c) It is enough to show that any subdigraph of $D\in D(n,p)$ on
$n^{1/9}$ vertices a.a.s.\  contains a cycle. To this end, let $C$ be
such a subdigraph. We can view $C$ as being sampled from
$D(n^{1/9},p)$. Using Theorem~\ref{Threshold}, we deduce that $p=n^{-1/9}$ is a
threshold for containing a cycle in $D(n^{1/9},p)$ (because the maximum
density of a cycle is 1), so because  $n^{-1/9}\log^2n\leq p=p(n)$, the
subdigraph $C$ a.a.s.\ contains a cycle. 

(d)  For an integer $k>n^{1/9}\log^2n/2$ and a set $S\subseteq V(D)$
with $|S|=k$,  let us enumerate $S$ as $\{1,2,\dots,k\}$. Let the
random variable $X$ count the number of arcs in the subdigraph induced
by $S$. Then $X=\sum_{1\leq i\neq j\leq k}X_{ij}$, where $X_{ij}$
counts the number of arcs (zero or one) from $i$ to $j$. Thus  
\begin{align*}
    E(X)=\sum_{1\leq i\neq j\leq k} E(X_{ij})=2\binom{k}{2}p.
\end{align*}
Using Corollary~\ref{Chernoffresult} with $\epsilon=1/\log n$, we have:
\begin{align*}
    P\Big(X\geq 2p\binom{k}{2}(1+1/\log n &)\; \text{or}\; X\leq 2p\binom{k}{2}(1-1/\log n)\Big)\\
    & \leq 2\exp{\left(-\frac{1}{3\log^2n}2\binom{k}{2}p\right)} \\ 
    & \leq 2\exp{\left(-\frac{1}{3\log^2n}k^2n^{-1/9}\log^2n\right)} \numberthis \label{ineq.1} \\
    & \leq 2\exp{\left(-\frac{k^2n^{-1/9}}{3}\right)}, 
    \numberthis \label{ineq.2}
\end{align*}
the estimate (\ref{ineq.1}) following from the hypothesis $p\geq n^{-1/9}\log^2n$. Now, suppose that $Y_t$ counts all the subsets of $V(D)$ of fixed size $t\geq k_0$ whose induced subdigraphs have at least $2p\binom{t}{2}(1+1/\log n)$ or at most $2p\binom{t}{2}(1-1/\log n)$ arcs. Then $Y=\sum_{t=k_0}^n Y_t$ counts all the subsets $U$ of size at least $k_0$ whose induced subdigraphs have at least $2p\binom{|U|}{2}(1+1/\log n)$ or at most $ 2p\binom{|U|}{2}(1-1/\log n)$ arcs. We have:
\begin{align*}
    E(Y) & =\sum_{t=k_0}^n E(Y_t)\\
    & \leq \sum_{t=k_0}^n 2\binom{n}{t}\exp{\left(-\frac{t^2n^{-1/9}}{3}\right)} \numberthis \label{ineq.3}
    \\ 
    & < \sum_{t=k_0}^n 2\Big(\frac{ne}{t}\Big)^t\exp{\left(-\frac{t^2n^{-1/9}}{3}\right)} \numberthis \label{ineqi.4}
    \\
    & =\sum_{t=k_0}^n 2\exp{\left(-t\log t+t\log n+t-\frac{t^2n^{-1/9}}{3}\right)}\\
    & = \sum_{t=k_0}^n 2 \exp{\Big(t\Big(\log n+1-\log t-\frac{tn^{-1/9}}{3}\Big)\Big)}\\
    &< 2\sum_{t=k_0}^n e^{-t}  \numberthis \label{ineqi.5}
    \\ & <2\sum_{t=k_0}^{\infty} e^{-t}=\frac{2e^{-k_0}}{1-e^{-1}}=o(1). 
    \numberthis \label{bound}
\end{align*}
The estimate (\ref{ineq.3}) follows from (\ref{ineq.2}), relation
(\ref{ineqi.4}) follows from the fact that
$\binom{n}{t}<(\frac{ne}{t})^t$, and (\ref{ineqi.5}) follows from 
the bound $\log n+1-\log t- \frac{tn^{-1/9}}{3}<-1$.
Using the bound (\ref{bound}) in Markov's inequality, we find that 
\begin{align*}
    P(Y=0)=1-P(Y\geq 1)\geq 1-E(Y)\rightarrow 1\; \text{as}\; n\rightarrow \infty.
 \end{align*}
So a.a.s.\ each set of $k\geq n^{1/9}\log^2n/2$ vertices induces a subdigraph with at most $2p\binom{k}{2}(1+1/\log n)=2p\binom{k}{2}(1+o(1))$ arcs.

(e) Let $S$ be a set of $k\geq k_1=n^{1/9}\log^2n$ disjoint pairs of vertices $\{v_i,w_i\}$, for $i\in [k]$ of $D\in D(n,p)$. Let $S'$ (the `contraction' of $S$) be the set obtained from $S$ by identifying $w_i$ with its corresponding $v_i$. For convenience, we enumerate $S'$ as $\{1,2,\dots,k\}$.
Now, suppose that $X$  counts the number of arcs (excluding loops and multiple arcs) in the subdigraph induced by $S'$. Then $X=\sum_{1\leq i\neq j\leq k}X_{ij}$, where $X_{ij}$  counts the number of arcs (zero or one) from $i$ to $j$ in the subdigraph induced by $S'$ (note that the sum is over ordered pairs). We have
\begin{align*}
    E(X_{ij})=P(X_{ij}=1)=1-P(X_{ij}=0)=1-(1-p)^4,
\end{align*}
so that 
\begin{align*}
    E(X)=\sum_{1\leq i\neq j\leq k} E(X_{ij})=2\binom{k}{2}\big[1-(1-p)^4\big].
\end{align*}
Using Corollary~\ref{Chernoffresult} with $\epsilon=1/\log n$, we have:
\begin{align*}
    P\Big[X\geq 2\binom{k}{2}\big(1-(1-p)^4\big)(1+1/\log n)\;\text{or}
    & \; X\leq 2\binom{k}{2}\big(1-(1-p)^4\big)(1-1/\log n)\Big] \\ 
    & \leq 2\exp{\left(-\frac{1}{3\log^2n}2\binom{k}{2}\big[1-(1-p)^4\big]\right)} \\
    & \leq 2\exp{\left(-\frac{1}{3\log^2n}2\binom{k}{2}p\right)} \numberthis \label{ineq.6}
    \\
    & \leq 2\exp{\left(-\frac{1}{3\log^2n}k^2n^{-1/9}\log^2n\right)} \\
    & = 2\exp{\left(-\frac{k^2n^{-1/9}}{3}\right)},
\end{align*}
where the estimate (\ref{ineq.6}) follows from the fact that $1-(1-p)^4\geq p$ for $0<p<1$.

Now, suppose that $Y_k$ counts all the sets with exactly $k$ disjoint pairs of vertices of $D$ whose contractions induce subdigraphs with at least $2\binom{k}{2}[1-(1-p)^4](1+1/\log n)$ or at most $2\binom{k}{2}[1-(1-p)^4](1-1/\log n)$ arcs (excluding loops and multiple arcs). Then $Y=\sum_{k=k_1}^nY_k$ counts all the sets with at least $k_1$ disjoint pairs whose contractions $U$  induce subdigraphs with at least $2\binom{|U|}{2}[1-(1-p)^4](1+1/\log n)$ or at most $2\binom{|U|}{2}[1-(1-p)^4](1-1/\log n)$ arcs. Arguing similarly to our estimates in part (d), we now have:
\begin{align*}
    E(Y) & =\sum_{k=k_1}^n E(Y_k)\\
    & \leq \sum_{k=k_1}^n 2\binom{n^2}{k}\exp{\left(-\frac{k^2n^{-1/9}}{3}\right)}\\ 
    & < \sum_{k=k_1}^n 2\Big(\frac{n^2e}{k}\Big)^k\exp{\left(-\frac{k^2n^{-1/9}}{3}\right)} \\
    & =\sum_{k=k_1}^n 2\exp{\left(-k\log k+2k\log n+k-\frac{k^2n^{-1/9}}{3}\right)}\\
    & = \sum_{k=k_1}^n  2\exp{\Big(k\Big(2\log n+1-\log k-\frac{kn^{-1/9}}{3}\Big)\Big)}\\
    & < 2\sum_{k=k_1}^n e^{-k}<2\sum_{k=k_1}^{\infty} e^{-k}=\frac{2e^{-k_1}}{1-e^{-1}}=o(1). \numberthis \label{bound2}
\end{align*}
Using the bound (\ref{bound2}) in Markov's inequality, we find that 

\begin{align*}
    P(Y=0)=1-P(Y\geq 1)\geq 1-E(Y)\rightarrow 1\; \text{as}\; n\rightarrow \infty.
 \end{align*}
So a.a.s.\ the contraction of each set $S$ of $k\geq n^{1/9}\log^2n$ disjoint pairs of vertices of $D$ induces a subdigraph with  $2\binom{k}{2}[1-(1-p)^4](1\pm 1/\log n)$ arcs (excluding loops and multiple arcs). It follows that in each set of $k$ disjoint pairs of vertices $\{\{v_i,w_i\}\}$, for $i\in \{1,2,\dots,k\}$ with $k\geq n^{1/9}\log^2n$,  there are  $2(1-(1-p)^4)\binom{k}{2}(1\pm o(1))$ pairs $(i,j)$ such that at least one of $v_iv_j,v_iw_j,w_iv_j,w_iw_j$ is an arc of $D$.
\end{proof}



\section{A.a.s.\ all digraphs are cores}
An  \textit{acyclic homomorphism} of a digraph $D$ to a digraph $C$, first defined in~\cite{Mark04},  is a function $\rho\colon V(D)\to  V(C)$ such that:
\begin{enumerate}[label=(\roman*)]
    \item for every arc $uv\in A(D)$, either $\rho(u)=\rho(v)$, or $\rho(u)\rho(v)$ is an arc of $C$; and 
    
    \item for every vertex $v\in V(C)$, the subdigraph of $D$ induced by $\rho^{-1}(v)$ is acyclic.
    
\end{enumerate}
For a more thorough treatment of graph and digraph homomorphisms, the reader is encouraged to consult~\cite{Hell04}.
We are now ready to state and prove the main result of this paper.
\begin{theorem}
\label{main-thm}
If $n^{-1/9}\log^2 n< p< 1-n^{-1/9}\log^2n$, and $D,C\in D(n,p)$, then a.a.s.\ every acyclic homomorphism $f\colon V(D)\to  V(C)$ is injective.
\label{Asymcore}
\end{theorem}
\begin{proof}
The bounds on $p$ imply that $D$ and $C$ a.a.s.\ satisfy properties
(a)--(e) in Lemma~\ref{Incomparable}. Suppose for a contradiction that
there exists an acyclic homomorphism $f\colon V(D)\to  V(C)$ that is
not injective. Then $f(x)=f(y)=z\in V(C)$ for some distinct vertices
$x,y\in V(D)$. Thus the set $A$ of vertices adjacent to either $x$ or
$y$ in $D$ must be mapped by $f$ to the set $B$ containing $z$ and
vertices adjacent to $z$. That is, if $A=N_D(x)\cup N_D(y)$ and
$B=N_C[z]$, then $f(A)\subseteq B$ (our notational convention being $N[z]=\{z\}\cup N(z)$).
Using (a) and (b) in Lemma~\ref{Incomparable}, a.a.s.\ we have
\begin{align*}
    |A|& \geq 2n(2p-p^2)(1-o(1))-np^2(2-p)^2(1+o(1))\\
    & \asymp \big(2np(2-p)-np^2(2-p)^2\big)(1-o(1))\\
    & = np(2-p)\big(2-p(2-p)\big)(1-o(1)),
 \end{align*}
 and 
  \begin{align*}
     & |f(A)|\leq |B|\leq n(2p-p^2)(1+o(1)).
\end{align*}
Thus a.a.s.\
\begin{align*}
    |A|-|f(A)|& \geq \big[np\big(2-p\big)\big(p^2-2p+2\big)\big](1-o(1))-np(2-p)(1+o(1))\\
    & \asymp \Big[np\big(2-p\big)\big(p^2-2p+1\big)\Big](1+o(1))\\
    & = np\big(2-p\big)\big(1-p\big)^2(1+o(1))\\
    & > \frac{1}{2}np(1-p)^2(1+o(1)) \\
    & \geq \frac{1}{2}n^{2/3}\log^6 n(1+o(1))  \numberthis  \label{bound3}\\ &
    \geq \frac{1}{2}n^{2/3}\log^2 n(1+o(1)), 
\end{align*}
where the bound (\ref{bound3}) follows from the fact that $p>n^{-1/9}\log^2n$ and  $1-p>n^{-1/9}\log^2 n$.
Because $f$ is an acyclic homomorphism, for any vertex $v\in V(C)$, the set $f^{-1}(v)$ is an acyclic set in $D$ so $|f^{-1}(v)|<n^{1/9}$  (part (c) of Lemma~\ref{Incomparable}).
Using the fact that $|A|-|f(A)|\geq n^{2/3}\log^2 n/2$ and $|f^{-1}(v)|<n^{1/9}$ shows that a.a.s.\ there are 
\begin{align*}
    k>\frac{|A|-|f(A)|}{n^{1/9}}>\frac{1}{2}n^{5/9}\log^2n>\frac{1}{2}n^{1/3}\log^2n>n^{1/9}\log^2n
\end{align*}
vertices $v_1,v_2,\dots,v_k\in f(A)$ such that $|f^{-1}(v_i)|\geq 2$. Using property (e) of Lemma~\ref{Incomparable}, we see that a.a.s.\ there are 
\begin{align*}
    2\big(1-(1-p)^4\binom{k}{2}\big)(1\pm o(1))
\end{align*}
arcs among the vertices in $\bigcup_{i=1}^k f^{-1}(v_i)\subseteq A$
and consequently among the vertices $v_1,v_2,\dots, v_k$. But part (d)
implies that there are at most $2p\binom{k}{2}(1+o(1))$ such
arcs. This gives our desired contradiction because
$2\big(1-(1-p)^4\binom{k}{2}\big)(1\pm o(1))>2p\binom{k}{2}(1+o(1))$.
\end{proof}

\begin{corollary}
If $n^{-1/9}\log^2 n< p< 1-n^{-1/9}\log^2n$, then a.a.s.\ a random digraph $D\in D(n,p)$ is a core.
\label{core}
\end{corollary}

\section{Acknowledgements}
 As noted on the title page, this work forms part of the first
 author's dissertation. He thanks his advisor, Mark Kayll, for
 his support, which came in various forms. 
 Thanks also to the referees for their contributions to improving the exposition.

 \bibliographystyle{siam}
\bibliography{Parsa}

\end{document}